\documentclass{article}

\usepackage{amsfonts}
\usepackage{amsmath}
\usepackage{amsthm}
\usepackage{amscd}
\usepackage{amssymb}
\usepackage[all]{xy}
\usepackage{geometry}
\xyoption{all}

\newtheorem{lemma}{Lemma}

\newtheorem{theorem}{Theorem}

\theoremstyle{remark}
\newtheorem{remark}{Remark}

\newtheorem{claim}{Claim}

\begin{document}

\title{The toric cobordisms}
\author{Alexandra MOZGOVA}
\date{}
\maketitle

\footnote{\noindent Laboratoire Emile Picard CNRS UMR 5580, Universit\'{e} Paul
Sabatier Toulouse III, 118, route de Narbonne, 31077 Toulouse,
France -- and -- Institute of Mathematics of Ukrainian National Academy of
Science, vul.~Tereschen\-kivska,~3, 252601 Kiev, Ukraine
\\
email : mozgova@picard.ups-tlse.fr
\\
This work was supported by French Government Grant \#19981314.
\\
Received March 1, 2001 and, in revised form,  August 23, 2002
\\
Communicated by : Ronald A. Fintushel}

\begin{abstract}
We introduce the notions of oriented and unoriented cobordisms in
the class of closed 3-manifolds fibered by tori $T^2$ and compute
the corresponding cobordism groups.
\end{abstract}

Two smooth compact $n$-manifolds, $M$ and $N$, are said to be
cobordant if their disjoint union $M\sqcup N$ is diffeomorphic to
the boundary of a smooth compact $(n+1)$-manifold $W$. Cobordism
is an equivalence relation on the set of smooth compact
$n$-manifolds, and the set of equivalence classes admits the
structure of an abelian group in which the operation is induced by
disjoint union. Other cobordism groups are obtained by restricting
the cobordism relation to classes of manifolds with additional
structure, such as an orientation, or a complex structure or a
spin structure. In this note we will consider the cobordism groups
of smooth compact $n$-manifolds which admit fibrations for which
the fiber is a $2$-dimensional torus.

An $n$-dimensional \textit{torus bundle} is a smooth fibration
$f:M^n\to B^{n-2}$ where the total space $M$ and the base $B$ are
smooth compact manifolds of dimension $n\ge2$ and $n-2$
respectively, and where the fibers of $f$ are diffeomorphic to the
$2$-dimensional torus $T^2$. By an \textit{oriented torus bundle} we
mean a torus bundle $f:M\to B$ together with a choice of an
orientation of $M$. If an (oriented) torus bundle $f:M\to B$ is
denoted by $M$, then $\partial M$ will denote the (oriented) torus
bundle obtained by restricting $f$ to the boundary of $M$ (with
the orientation induced from that of $M$); we will use $M^-$
to denote the oriented torus bundle consisting of the fibration
$f:M\to B$ together with the opposite choice of orientation of
$M$. If $M_1$ and $M_2$ are (oriented) torus bundles, then we will
write $M_1=M_2$ if there is a fiber-preserving
(orientation-preserving) diffeomorphism from $M_1$ to $M_2$.

Two $n$-dimensional torus bundles $M_1$ and $M_2$ are said to be
\textit{toric cobordant} if there exists an $(n+1)$-dimensional torus
bundle $W$ such that $\partial W = M_1\sqcup M_2$. If $M_1$ and
$M_2$ are oriented torus bundles, then we will say that they are
\textit{oriented toric cobordant} if there exists an oriented torus
bundle $W$ such that $\partial W = M_1\sqcup M_2^-$. In either
case $W$ will be called a \textit{cobordism} from $M_1$ to $M_2$.
Toric cobordism and oriented toric cobordism are equivalence
relations. We will use the notation $M_1 \sim M_2$ to indicate
that the (oriented) torus bundles $M_1$ and $M_2$ are (oriented)
toric cobordant, and in either case the equivalence class of $M$
will be denoted $M^\omega$. The equivalence classes of
$n$-dimensional torus bundles under toric cobordism form an
abelian group $\Omega_n^{T^2,unor}$ under the operation defined by
$M_1^\omega + M_2^\omega = (M_1\sqcup M_2)^\omega$. The
equivalence classes of $n$-dimensional oriented torus bundles
under oriented toric cobordism form an abelian group
$\Omega_n^{T^2,or}$ under the operation defined by $M_1^\omega +
M_2^\omega = (M_1\sqcup M_2^-)^\omega$. 
The 
unoriented toric cobordism groups of $n$-manifolds are isomorphic
to the unoriented bordism groups $\Omega_{n-2}(B
GL(2,\mathbb{Z}))$ of maps from $B^{n-2}$ to $B Diff(T^2)$.  Lemma 1 below
shows that the oriented toric cobordism  groups of 3-manifolds are isomorphic
to the oriented bordism groups $\Omega_1(BGL(2,\mathbb{Z}))$ (pointed out by
A.~Marin).

The rest of this note will be devoted to the proof of the
following:
\begin{theorem}\label{main}
The third group of oriented toric cobordisms $\Omega_{3}^{T^2,or}$
is isomorphic to $\mathbb{Z}_{12}$ and is generated by the class
of the torus bundle with the monodromy
\begin{math} \bigl(
   \begin{smallmatrix}
   1 & 1 \\
   0 & 1
  \end{smallmatrix}
  \bigr)
\end{math}.
The third group of unoriented toric cobordisms
$\Omega_{3}^{T^2,unor}$ is isomorphic to $\mathbb{Z}_{2} \oplus
\mathbb{Z}_{2}$ and is generated by the classes of the torus
bundles with the monodromies
\begin{math} \bigl(
   \begin{smallmatrix}
   0 & -1 \\
   1 & 0
  \end{smallmatrix}
  \bigr)
  \end{math} and
 \begin{math} \bigl(
   \begin{smallmatrix}
   0 & 1 \\
   1 & 0
  \end{smallmatrix}
  \bigr)
\end{math}.
\end{theorem}

\section*{The toric cobordisms of $3$-manifolds}
The quotient $M_{\varphi}=T^{2} \times I /(x,0) \sim
(\varphi(x),1)$ where $\varphi$ is a diffeomorphism of $T^{2}$ is
a torus bundle over the circle with monodromy $\varphi$. If the
diffeomorphisms $\varphi_1, \varphi_2$ are isotopic, then
$M_{\varphi_1} = M_{\varphi_2}$. A diffeomorphism of $T^2$ is
determined up to isotopy by its induced map on the first integral
homology group, and hence the diffeotopy group is isomorphic to
$GL(2, \mathbb{Z})$.

Let $G$ denote the group $GL(2,\mathbb{Z}$); let $G'$ be its
commutator subgroup and $G^2$ be the subgroup of $G$ generated by
the squares of the elements of $G$.
\begin{claim}\label{t-or}
Let $M_\varphi$ be a 3-dimensional torus bundle defined as above.
There exists a $4$-dimensional manifold $W$ fibered by tori over
an orientable surface with $\partial W = M_\varphi$ if and only if
$\varphi \in G'$. Such a $W$ is orientable if and only if
$\varphi$ can be written as $\varphi =
\prod_{j=1}^{g}[\varphi_{2j- 1},\varphi_{2j}]$ where $\operatorname{det}
\varphi_i = 1$ for $i=1, \ldots ,2g$.
\end{claim}
\begin{proof}
The $T^2$-bundles over a finite cell complex $X$ are classified by
homotopy class of maps from $X$ to $B Diff(T^2)$. Such a class
determines a conjugacy class of homomorphisms \[\pi_1(X) \to
\pi_1(B Diff(T^2)) \cong \pi_0(Diff(T^2)) \cong \ G.\] In our case
$X$ is a surface with non-empty boundary, so $\pi_1(X)$ has
cohomological dimension 1 (being free). Thus $T^2$-bundles over a
surface with non-empty boundary are in bijection with the
homomorphisms $\pi_1(X) \to G$.

The manifold $W$ we are looking for exists if and only if there is
a commutative diagram \[\xymatrix{ \pi_1(S^1) \cong \mathbb{Z}
\ar[r] \ar[d]_{i_{\ast}} &  G
\\ \pi_1(B^2) \cong F_r \ar[ur] & } \]
where $i_{\ast}$ is induced by inclusion of the boundary, i.e. $Im
\ i_{\ast} \subset (F_r)'$, hence we have shown the first part of the
claim.

To see when $W$ is orientable, use the standard construction of an
oriented surface from the disk $D^2$, by identifying some
$1$-disks $I_i$ on its boundary. Then we can explicitly construct
$W$ starting from $D^2 \times T^2$ by gluing $(I_i \times T^2)$'s
on its boundary. As $D^2 \times T^2$ is oriented and the gluing
must create no orientation-reversing loop, it is not hard to see
that the condition $\operatorname{det}\varphi_i=1$ assures the
orientability of $W$.
\end{proof}
Reasoning similarly, we get for fiber bundles over a
non-orientable surface
\begin{claim}\label{t-nonor}
Let $M_\varphi$ be a 3-dimensional torus bundle. There exists a
$4$-dimensional manifold $W$ fibered by tori over a non-orientable
surface with $\partial W = M_\varphi$ if and only if $\varphi \in
G^2$. Such a $W$ is orientable if and only if $\varphi$ can be
written as $\varphi= \prod_{j=1}^{k} \varphi^2_{j}$ where $det\
{\varphi}_i=-1$ for $j=1, \ldots ,k$.
\end{claim}

\begin{claim}\label{many-or}
Let $\varphi_{1}\dots\varphi_{n} \in G=GL(2,\mathbb{Z})$ and
$M_{\varphi_{1}}, \ldots, M_{\varphi_{n}}$ be the corresponding
3-dimensional torus bundles. There exists a $4$-dimensional
manifold $W$ fibered by tori over an orientable surface with
$\partial W = M_{\varphi_{1}} \sqcup \ldots \sqcup
 M_{\varphi_{n}}$ if and only if $\prod_{i=1}^{n} \varphi_{i} \in
G'.$
\end{claim}
\begin{proof}
The proof is immediately obtained from Claim \ref{t-or} and the
fact that the disjoint union $M_{\varphi_{1}} \sqcup \ldots
\sqcup M_{\varphi_{n}}$ is cobordant to $M_{\psi}$ with
$\psi=\prod_{i=1}^{n} \varphi_{i},$ by a toric cobordism with base
the sphere $S^2$ with $n+1$ holes (it can be constructed similarly
to the proof of Claim \ref{t-or}).
\end{proof}
\begin{claim}\label{many-nonor}
For $\varphi_{1}\dots\varphi_{n} \in G$ and $M_{\varphi_1},
\ldots, M_{\varphi_k}$ as above, there exists a $4$-dimensional
manifold $W$ fibered by tori over a non-orientable surface with
boundary $\partial W = M_{\varphi_{1}} \sqcup \ldots \sqcup
 M_{\varphi_{n}}$ if and only if $\prod_{i=1}^{n} \varphi_{i} \in
G^2.$
\end{claim}

\begin{lemma}\label{or=non}
There exists an oriented toric cobordism with an orientable base
between $M_{\varphi}$ and $M_{\psi}$ if and only if there also
exists an oriented toric cobordism with a non-orientable base
between them.
\end{lemma}
\begin{proof}
By $G_{-}^2$ we denote the subgroup of $GL(2,\mathbb{Z})$
generated by the squares of matrices with negative determinant 
\[G_{-}^2=\langle \{ a_1^2a_2^2 \dots a_k^2 \mid det\ a_i=-1 \}
\rangle.\] It is evident that $G_{-}^2 \subset SL(2,\mathbb{Z})$
and is normal in it. We show that $G_{-}^2 = (SL(2,\mathbb{Z}))'$
and this implies the claim. We use the following presentations of
$GL(2,\mathbb{Z})$ and $SL(2,\mathbb{Z})$; see (\cite{Z2}, 2.23).
For $A=$
\begin{math} \bigl(
   \begin{smallmatrix}
   0 & -1 \\
   1 & 0
  \end{smallmatrix}
  \bigr)
  \end{math},
 $B=$
\begin{math} \bigl(
   \begin{smallmatrix}
   0 & 1 \\
   -1 & 1
  \end{smallmatrix}
  \bigr)
\end{math} and $R=$
\begin{math} \bigl(
   \begin{smallmatrix}
   0 & 1 \\
   1 & 0
  \end{smallmatrix}
  \bigr)
\end{math},
\[GL(2,\mathbb{Z})=\langle A,B,R | A^2=B^3,
A^4=R^2={(RA)}^2={(RB)}^2=1 \rangle,\] \[SL(2,\mathbb{Z})=\langle
A,B | A^2=B^3, A^4=1 \rangle.\] The commutator subgroup
$(SL(2,\mathbb{Z}))'$ of $SL(2,\mathbb{Z})$ is a free group of
rank 2 generated by $\lbrack A,B \rbrack = {(ARB^{-1})}^2 \in
G_{-}^2$ and $\lbrack A,B^{-1} \rbrack =(B^{-1}RB)^2 \in G_{-}^2.$
Thus $(SL(2,\mathbb{Z}))' \subset G^2_{-}$, and as
$(SL(2,\mathbb{Z}))' \vartriangleleft SL(2,\mathbb{Z})$, we have
$(SL(2,\mathbb{Z}))' \vartriangleleft G^2_{-}$.

By using the relations $RA=A^{-1}R$ and $RB=B^{-1}R$, each element
$a \in GL(2,\mathbb{Z})$ can be written in the normal form
$R^{\varepsilon}A^{k_1}B^{l_1} \dots A^{k_n}B^{l_n}$ where
$\varepsilon \in \{0;1\}$. If $det \ a=-1$, the element $a$ can be
written in the form $RA^{k_1}B^{l_1} \dots A^{k_n}B^{l_n}$.

Now $G^2_{-} / (SL(2,\mathbb{Z}))'= \langle \{ (RA^iB^j)^2 \ |\
A^2=B^3,(RA)^2=(RB)^2=A^4=R^2=1, AB=BA \}\rangle =1.$ Thus,
$G^2_{-}=(SL(2,\mathbb{Z}))'$.
\end{proof}

\begin{proof}[Proof of Theorem \ref{main}] It follows from
Lemma \ref{or=non} that for $\varphi, \psi \in SL(2,\mathbb{Z})$
there exists an oriented toric cobordism between the torus bundles
$M_{\varphi}$ and $M_{\psi}$ if and 
only if $\varphi {\psi}^{-1}
\in (SL(2,\mathbb{Z}))'.$ Thus \[\Omega_{3}^{T^2,or} \cong
SL(2,\mathbb{Z}) / (SL(2,\mathbb{Z}))'= \langle A,B | A^2=B^3,
A^4=1, AB=BA \rangle \cong \mathbb{Z}_{12}.\] The generator of
$\Omega_{3}^{T^2,or}$ is the conjugacy class of the element
$B^{-1}A^{-1}=$
\begin{math} \bigl(
   \begin{smallmatrix}
   1 & 1 \\
   0 & 1
  \end{smallmatrix}
  \bigr)
\end{math}.

The subgroup $G'$ lies in the subgroup $G^2$. Thus, if there
exists an unoriented toric cobordism with an orientable base
between $M_{\varphi}$ and $M_{\psi}$, then there also exists an
unoriented toric cobordism with a non-orientable base between
them. Hence, Claim \ref{many-nonor} lets us calculate the
third group of unoriented toric cobordims,
$\Omega_{3}^{T^{2},unor}$, as well. Thus, \[\Omega_{3}^{T^2,unor}
\cong G / G^2= \langle A,B,R | A^2=B^3,
A^4=R^2={(RA)}^2={(RB)}^2=1, \] \[ AB=BA, AR=RA, BR=RB, A^2=B^2=1
\rangle \cong \mathbb{Z}_{2} \oplus \mathbb{Z}_{2}.\] The
generators here are the conjugacy classes of $A=$
\begin{math} \bigl(
   \begin{smallmatrix}
   0 & -1 \\
   1 & 0
  \end{smallmatrix}
  \bigr)
  \end{math} and
 $R=$
\begin{math} \bigl(
   \begin{smallmatrix}
   0 & 1 \\
   1 & 0
  \end{smallmatrix}
  \bigr).
\end{math}
\end{proof}

\begin{remark}
The toric cobordism is a cobordism of manifolds with
some fixed torus bundle structures. In some cases, for
$M_{\varphi},\ M_{\psi}$ that are not oriented toric cobordant, we can
choose other torus bundle structures on their total spaces in such
a way that they become oriented toric cobordant. For example, take a
$\varphi\in SL(2,\mathbb{Z})$ 
such that $\varphi^2 \nsubseteq  (SL(2,\mathbb{Z}))'$.  
Then, the corresponding
$M_{\varphi}$ is not oriented toric cobordant to $M_{\varphi}^{-}$, but
is oriented toric cobordant to $M_{\gamma{\varphi}{\gamma}^{-1}}$ for $\gamma$
with $\det \gamma=-1$.  The reason for this is that the total spaces of $M_
\varphi$ and $M_{\gamma\varphi\gamma^{-1}}$ are homeomorphic by a fiber
preserving homeomorphism inducing the orientation reversing map of the basis
$S^1$ (see \cite{BZ} or \cite{Hatcher}).
\end{remark}

\begin{remark} Let $\varphi \in SL(2, \mathbb{Z})$ and $W^4$ be a
toric cobordism between $M_{\varphi}$ and $\varnothing$ that has
an orientable base of genus $g$. Then $\varphi$ is a product of
$g$ commutators. By taking $\psi', \psi''$ such that $\lbrack
\psi', \psi'' \rbrack =1$, one can write $\varphi$ as a product of
$g+1$ commutators and so one can construct another toric cobordism
between $M_{\varphi}$ and $\varnothing$ with orientable base of
genus $g+1$. Thus, quite naturally we come to the following
question: what is the minimal genus of the orientable base of
$W^4$?

For this, we can utilize Culler's algorithm (\cite{CU}) which
determines, for finite groups $A$ and $B$ and $a \in (A \ast B)'$,
the minimal number of elements of $A \ast B$ required to represent
$a$ as a product of their commutators. In order to extend Culler's
algorithm from free products to $SL(2, \mathbb{Z})={\mathbb
Z}_4\ast_{{\mathbb Z}_2}{\mathbb Z}_6$, consider the projection
homomorphism $\alpha : {\mathbb Z}_4 \ast_{{\mathbb Z}_2} {\mathbb
Z}_6 \to {\mathbb Z}_2 \ast {\mathbb Z}_3.$ As Culler's algorithm
can be applied to the group ${\mathbb Z}_2 \ast {\mathbb Z}_3$, it
remains to note that the restriction of the homomorphism $\alpha$
to the commutator subgroup
\[\alpha'= \alpha |_{({\mathbb Z}_4
\ast_{{\mathbb Z}_2} {\mathbb Z}_6)'}: ({\mathbb Z}_4
\ast_{{\mathbb Z}_2} {\mathbb Z}_6)' \longrightarrow ({\mathbb
Z}_2 \ast {\mathbb Z}_3)'\] is an isomorphism.
\end{remark}

\section*{Acknowledgment}
The author is grateful to Gilbert Levitt for pointing out Culler's
paper. It is also a pleasure to thank the referee for constructive
comments.

\bibliographystyle{amsplain}

\end{document}